\documentclass[12pt, a4paper]{article}
\usepackage{mathtext}
\usepackage{cmap}					
\usepackage[T2A]{fontenc}			
\usepackage[utf8x]{inputenc}			
\usepackage[english]{babel}	
\usepackage{indentfirst}
\usepackage{amsfonts}
\usepackage{xcolor}
\usepackage{amsthm}
\usepackage{amsmath}
\usepackage{blindtext}
\usepackage{graphicx}
\usepackage{indentfirst}
\usepackage{amscd,amssymb}
\usepackage{amsfonts,amsmath,array}
\usepackage{mathtools}
\usepackage{url}

\newtheorem{theorem}{Theorem}[section]
\newtheorem{lemma}[theorem]{Lemma}
\newtheorem{remark}[theorem]{Remark}
\newtheorem{example}[theorem]{Example}
\newtheorem{definition}[theorem]{Definition}
\newtheorem{proposal}[theorem]{Proposition}
\newtheorem{consequence}[theorem]{Сorollary}
\usepackage{titlesec}
\usepackage{listings}

\usepackage{setspace}
\usepackage[12pt]{extsizes}
\usepackage[left = 30mm, top=20mm, right = 10mm, bottom=20mm, nohead, nofoot]{geometry}

\usepackage{fancyhdr}
\title{On the Theorem of Gauss--Lucas for quaternions}
\author{ I. Emizh$^{1,2}$, A. Guterman$^{3}$  \\
        \small $^{1}$  Caucasus Mathematical Center of Adyghe State University, Maykop 385000, Russia \\
        \small  \small $^{2}$ Moscow Institute of Physics and Technology, Dolgoprudny 141700, Russia \\
        \small $^{3}$  Bar-Ilan University, Ramat-Gan 5290002, Israel}
\date{}

\begin{document}
\newpage
\maketitle

    \begin{abstract}
        It is proved that the roots of the derivative of a polynomial with quaternionic coefficients belong to the union of the intersections of   sets defined in terms of certain projections of a polynomial. 
        The result strengthens the quaternion version of Gauss-Lucas  theorem, proved by Ghiloni
        and Perotti in 2018.  


    \end{abstract}

\section{Introduction}
    The famous Gauss -- Lucas theorem describes the geometric relationship between the roots of a complex polynomial and the roots of its derivative.
    
    \begin{theorem}[Gauss -- Lucas]\label{GL}
        Let $f \in \mathbb{C}[x]$ be an arbitrary polynomial of degree no less than 2 of one variable with complex coefficients. Then the roots of the derivative $f'$ of the polynomial $f$ belong to the convex hull of the roots of the original polynomial, see~\cite[\S\:2.1]{PV}, here $\mathbb{C}$ denote the field of complex numbers. 
    \end{theorem}

    Various results on functions of a complex variable have been successfully generalized to functions on quaternions for about a hundred years. 
    For example, in the works \cite{RF, RFS} analogues of the Cauchy-Riemann equations, the Cauchy integral formula, and the Laurent expansion have been proposed.  The fundamental theorem of algebra for polynomials over quaternions was proved in paper \cite{SEIN}. A number of works have been devoted to the study of the properties of the roots of polynomials over quaternions and other division rings, see \cite{ACAG, GBMT, VF, NI} and the references therein. One of the important and promising results for generalization is the Gauss-Lucas theorem. This theorem provides a geometric relationship between the roots of a polynomial and the roots of its derivative, which in particular offers methods for estimating the modulus of polynomial roots, thereby opening up opportunities for numerical methods. Unsuccessful attempts to generalize it for quaternions have been made in the works \cite{VF, SGJOIS}. However, in the work of Ghiloni and Perotti \cite{GRPA} an example was constructed showing that a direct generalization for an arbitrary polynomial is impossible. Moreover, a class of polynomials was distinguished for which a direct generalization is valid. In particular, Ghiloni and Perotti  introduced a new set that depends on the roots of the polynomial, which contains the roots of the derivative of this polynomial in the general case. An estimate of the absolute values of the roots of the derivative of the polynomial in terms of the roots of the original polynomial was obtained, see \cite[Proposition~6]{GRPA}.

    The main purpose of this work is to strengthen the results from \cite{GRPA}. It is proved that the roots of the derivative of a polynomial with quaternionic coefficients belong to the union of the intersections of convex hulls of the roots of some projections of this polynomial. It is shown that the set considered in this paper is contained in the set constructed in the paper \cite{GRPA}. Moreover examples of polynomials for which this inclusion is strict are presented. We also obtained an enhancement of the estimation of the absolute value of the roots of the derivative of a polynomial through the roots of the original polynomial proved in~\cite{GRPA}.    
    
    We are going to introduce the basic concepts used in our paper  following the notations proposed in~\cite{GRPA}.
    Similar to the case of polynomials over $\mathbb{R}$ or over $\mathbb{C}$, polynomials over $\mathbb{H}$ can be considered not only as formal expressions, but also as functions. That is, the polynomial $P(x) = \sum\limits_{t = 0}^{n}a_{t}x^{t} \in \mathbb{H}[x]$ is   a function $P(x): \mathbb{H} \to \mathbb{H}$, such that for each element $h \in \mathbb{H}$ its image is
    \begin{equation*} P(h) = \sum\limits_{t = 0}^{n}a_{t}h^{t}.\end{equation*}
    Note that $a_t$ in general does not commute with $h$. 
    Denote by $R(P)$ the set of roots of a polynomial $P$, that is the set of quaternions $h \in \mathbb{H}$, such that  $P(h) = 0$. 

    Consider a polynomial $P(x) \in \mathbb{H}[x]$ and an imaginary quaternion $I \in \mathcal{S}$ of the unit norm, where $\mathcal{S}$ is the imaginary unit sphere. Let $\mathbb{R}[I] \subseteq \mathbb{H}$ be the plane generated by $1$ and $I$. In second section we will show that $\mathbb{R}[I]$ can be considered as a complex plane, see Proposition \ref{local complex}. For $P(x)$ denote by $P_{I}(x)$ the projection of this a polynomial onto $\mathbb{R}[I]$, see Definition \ref{def poly proj}. In the third section, it will be shown that the set of such projections can be considered as a ring $\mathbb{C}[x]$. Denote by  $\mathbf{conv} (E)$ convex hull of the set $E \subseteq \mathbb{H}$. 
    
    The following result will be in the center of our considerations.
    \begin{theorem}\label{QGL}\cite[Theorem 1]{GRPA}
        For any polynomial $P \in \mathbb{H}[x]$ with degree no less than 2
        $$R(P') \subseteq  \mathfrak{sn}(P).$$
    \end{theorem}    
    
    Here $\mathfrak{sn}(P) = \bigcup\limits_{I \in \mathcal{S}} \mathbf{conv} (R(P_{I})) $ this is a   convex hull proposed in~\cite{GRPA} and called there a snail.  
    
    In order to formulate our result  we  introduce another projection $P_{I}^{\perp}(x)$ for the polynomial $P(x)$, see Definition~\ref{def poly proj}. The main result of our work is the following:
    
    \begin{theorem}\label{SQGL}
        For any polynomial $P \in \mathbb{H}[x]$ it holds that
        $$R(P') \subseteq \bigcup\limits_{I \in \mathcal{S}} \mathbf{conv} (R(P_{I})) \cap \mathbf{conv} (R(P_{I}^{\perp})). $$  
    \end{theorem}
    
    Our paper is organized as follows. In Section 2 we provide necessary preliminary information about quaternions.  Section 3 is devoted to the theory of polynomials with quaternionic coefficients. In   Section 4 the proof of Theorem \ref{SQGL} is provided. 
    In   Section 5 we apply this  theorem to estimate the absolute value of the roots of the derivative of a polynomial and show that Theorem \ref{SQGL} indeed  strengthens the  estimate from \cite[Proposition~6]{GRPA}.

\section{Preliminary information about quaternions}

    Quaternions form an algebra over $\mathbb{R}$, denoted by $\mathbb{H}$ and generated by the real unit  $1$ and three imaginary units $ i,j,k,$   which satisfy the relations:
    $$i^2 = j^2 = k^2 = ijk = -1.$$
    
    Let us define quaternion as $h = \alpha_0 + \alpha_1 i + \alpha_2 j + \alpha_3 k$. The real and imaginary parts of the quaternion $h$ are defined as $Re(h)  = \alpha_0$ and $Im(h) = \alpha_1 i + \alpha_2 j + \alpha_3 k$, respectively. For quaternion $h$ we define its conjugate by $\overline{h} = \alpha_0 - \alpha_1 i - \alpha_2 j - \alpha_3 k = Re(h) - Im(h)$.
    
    On the algebra of quaternions a positive definite quadratic form $Q(h) = \|h\| = h \cdot \overline{h} = \alpha_0^2 + \alpha_1^2 + \alpha_2^2 + \alpha_3^2 \ge 0$ is well defined and it determines  the inner product space structure in the usual way. Namely, for two  quaternions 
    $p = \beta_0 + \beta_1 i + \beta_2 j + \beta_3 k, \,\,\, q = \gamma_0 + \gamma_1 i + \gamma_2 j + \gamma_3 k\in \mathbb{H}$
    their scalar product is defined to be 
    \begin{equation} (p,q) = \beta_0 \cdot \gamma_0 + \beta_1 \cdot \gamma_1  + \beta_2 \cdot \gamma_2 + \beta_3 \cdot \gamma_3.\label{2} \end{equation}
    The orthogonality relation in our paper is understood in the sense of this scalar product. For any nonzero quaternion $h$ there is the unique multiplicative inverse element which can be computed as   $h^{-1} = \frac{\overline{h}}{\|h\|}$.
    
    Let  $\mathcal{S} = \{ h \in \mathbb{H} \,\,\, | \,\,\, \|h\| = 1, \,\,\, Re(h) = 0\}$ be the imaginary unit sphere and $I \in \mathcal{S}$ be an arbitrary fixed quaternion lying on this sphere. Denote by $$\mathbb{R}[I] = \{a + Ib \,\,\, | \,\,\, a,b \in \mathbb{R} \} \subseteq \mathbb{H}$$ 
    the plane generated by $1$ and $I$. This plane is a subalgebra of the algebra $\mathbb{H}$. The following result illustrates  the structure of the subalgebra $\mathbb{R}[I]$ see \cite{MBPS}.

    \begin{proposal}\label{local complex} \cite[Theorem 4.7]{MBPS}
        For any quaternion  $I \in \mathcal{S}$ the algebra $\mathbb{R}[I] \subseteq \mathbb{H}$ generated by $1$ and $I$ is isomorphic to the algebra~$\mathbb{C}$.
    \end{proposal}

    We will need the following remark for further considerations. 
     \begin{remark}\label{complex plane} 
    The isomorphism $\phi : \mathbb{R}[I] \to \mathbb{C}$ can be explicitly defined as follows: for the element $a + Ib \in \mathbb{R}[I]$, where $a,b \in \mathbb{R}$ put $\phi (a + Ib) = a + ib$. A direct computations show that this is indeed an isomorphism of algebras, moreover, this isomorphism preserves the norm.
    \end{remark}

    Proposition \ref{local complex} shows that any element $I \in \mathcal{S}$ can be considered as an imaginary unit and the plane $\mathbb{R}[I]$ is indeed the complex plane.

    We decompose the quaternion into a sum of the form $z_{1} + Jz_{2}$, where $J$ a fixed quaternion of the unit norm, and quaternions $z_{1}, z_{2}$ lie in the plane $\mathbb{R}[I]$. We require that $z_{1}$ is an orthogonal projection of $z$ on $\mathbb{R}[I]$ and $Jz_{2}$ is a projection of $z$ onto an orthogonal complement to $\mathbb{R}[I]$. 
    
    \begin{lemma}\label{rem 3 I}
        For a given quaternion $I \in \mathcal{S}$ there is a quaternion $J \in \mathcal{S}$, such that for any $h \in \mathbb{H}$ there is a unique decomposition of the form $h = z_{1} + Jz_{2}$, where $z_1, z_2 \in \mathbb{R}[I]$. Moreover, the orthogonal projection of $h$ on $\mathbb{R}[I]$ is equal to $z_1$ and the orthogonal complement is $Jz_2$. 
    \end{lemma}
    
    Proof:
    By the definition, each quaternion can be represented as $h = (\alpha_{0} + \alpha_{1}i) + j(\alpha_{2} + \alpha_{3}i)$, where $\alpha_{0}, \alpha_{1}, \alpha_{2}, \alpha_{3} \in \mathbb{R}$. For the selected $I \in {\cal{S}}$ there is  an element $g \in \mathbb{H}^{*}$ such that $I = gig^{-1}$, since  $Re(I) = 0$ and $\|I\| = 1$, see \cite[Proposition 3.6]{ACAG}. Conjugation by an element $g$ is an automorphism of $\mathbb{H}$. So, there is a unique set of elements $\beta_{0}, \beta_{1}, \beta_{2}, \beta_{3} \in \mathbb{R}$, satisfying the equality: $h = (\beta_{0} + \beta_{1} I) + J(\beta_{2} + \beta_{3}I)$, where $J = gjg^{-1}$. Moreover, the automorphism obtained by conjugation on $g$ is a rotation of the imaginary part \cite[\S 3.1]{CJSD}, therefore, conjugation on g preserves the scalar product. More strictly, we use two properties that are proved by direct calculation
    \begin{enumerate}
        \item $(h_{1}, h_{2}) = \frac{1}{2}(\|h_{1} + h_{2}\| - \|h_{1}\| - \|h_{2}\|),$
        \item $\|ghg^{-1}\| = g\|h\|g^{-1}.$
    \end{enumerate}
    Then we have
    $$(gh_{1}g^{-1}, gh_{2}g^{-1}) = \frac{1}{2}(\|g(h_{1} + h_{2})g^{-1}\| - \|gh_{1}g^{-1}\| - \|gh_{2}g^{-1}\|) = $$
    $$= g(\frac{1}{2}(\|h_{1} + h_{2}\| - \|h_{1}\| -\|h_2\|))g^{-1} = g(h_{1}, h_{2})g^{-1} = (h_{1}, h_{2}).$$
    It follows  that this automorphism transfers $\mathbb{R}[i]$ to $\mathbb{R}[I]$. Also its orthogonal complement, denoted by $\mathbb{R}[i]^{\perp}$, which is a plane generated by $j$ and $k$, is mapped into the  orthogonal complement $\mathbb{R}[I]^{\perp}$ which is a plane generated by $J$ and $K = gkg^{-1}$. Thus, the projection $h$ on $\mathbb{R}[I]$ is $\beta_{0} + \beta_{1} I$ and the projection on it's orthogonal complement is equal to $J(\beta_{2} + \beta_{3}I)$, respectively.

    \qed
    
    Let us illustrate this lemma by considering the following example
    \begin{example}
        Given a quaternion $h = 1 + 11i -2j -k$, it can be represented as 
        $$h = (1 + 11i) + j(-2 + i).$$   
        Fix $ I = \frac{1}{3}i + \frac{2}{3}j - \frac{2}{3}k$. Let us find an element $g$, such that $I = gig^{-1}$. Then we calculate the elements $g = 1 + i + j$ and $g^{-1} = \frac{1 - i - j}{3}$. Therefore $J = gjg^{-1} = \frac{2}{3}i + \frac{1}{3}j + \frac{2}{3}k$. Let's find the coefficients  $\beta_{0}, \beta_{1}, \beta_{2}, \beta_{3}$. To calculate the coefficients, we use the fact that 1, I, J, IJ form an orthonormal basis:
        $$\beta_{0} = (1, h) = 1$$
        $$\beta_{1} = (I, h) = 3$$
        $$\beta_{2} = (J, h) = 6$$
        $$\beta_{3} = (JI, h) = -9$$
        Finally $$h = (1 + 3I) + J(6 - 9I).$$
        Considering that $1$ and $I$ are orthogonal quaternions with the unit norm, we obtain that the orthogonal projection $h$ on $\mathbb{R}[I]$ is
        $$(1,h) + (I,h) \cdot I = 1 + 3I.$$
        Accordingly, orthogonal quaternions of the unit norm $J$ and $JI = -\frac{2}{3}i + 
        \frac{2}{3}j + \frac{1}{3}k$ generate the plane $\mathbb{R}[I]^{\perp}$. Then the orthogonal complement along $\mathbb{R}[I]$ is equal to
        $$(J,h) \cdot J + (JI,h) \cdot JI = J(6 - 9I).$$
    \end{example}

    \begin{remark}
        Observe that for the quaternion $h$ from the previous example we have simultaneously
        $$h = (1 + 3I) + J(6 - 9I) = (1 + 3I) - J(-6 + 9I).$$
        We remark that Lemma \ref{rem 3 I} is not stating the uniqueness of the choice of such a quaternion $J$.
    \end{remark}
    
\section{Preliminary information about polynomials}

    \begin{definition}\label{defolt def}
        \begin{enumerate} 
            \item Polynomial over quaternions is the formal expression of the form 
            $$f(x) = a_{n} x^n + a_{n - 1} x^{n - 1} + \ldots + a_{0},$$
            where $a_{t} \in \mathbb{H}, \,\,\, t = \overline{0, n}$. The set of polynomials will be denoted by $\mathbb{H}[x]$.

            \item The derivative of a polynomial $f \in \mathbb{H}[x]$ is the expression of the form 
            $$f'(x) = n \cdot a_{n} x^{n-1} + (n - 1) \cdot a_{n - 1} x^{n - 2} + \ldots + a_{1}.$$

            \item The root of the polynomial $f \in \mathbb{H}[x]$ is   a quaternion $h$ such that   
            $$f(h) = a_{n} h^n + a_{n - 1} h^{n - 1} + \ldots + a_{0} = 0.$$
            The set of roots of $f$ will be denoted by $R(f)$.
        \end{enumerate}
    \end{definition}
    Let us illustrate each of these definitions with an example.
    \begin{example} \label{primer IR}
        \begin{enumerate}
            \item    $f(x) = x^3 + i x^2 + j x + k$ -- a polynomial over quaternions.

            \item $f'(x) = 3x^2 + 2i x + j$ --   the derivative of the polynomial $f(x)$ from the previous paragraph.
            
            \item  Let $p(x)=x^2 + 1$. Then $R(p(x)) = \{ h \in \mathbb{H} \,\,\, | \,\,\, \|h\|=1, \,\,\, Re(h) = 0\}$.
        \end{enumerate}
    \end{example}
    
     In the Example \ref{primer IR}.3 the set of different roots of the polynomial $p(x)$ is infinite. Moreover, if we decompose $p(x) = (x - j)(x + j)$, then the roots of $p(x)$ will not be the union of the roots of linear factors. Since, for example, $i$ which is the root of the polynomial $p(x)$ is neither the root of $x - j$ nor the root of $x + j$. This is due to the fact that if we consider $x$ as an unknown quaternion, and not as a formal symbol, then $(x - j)(x + j) = x^2 - jx + xj + 1$ not equal $x^2 + 1$ due to the lack of commutativity. Therefore, if we want to extract the roots of a product from its multipliers, then we need to follow the order of multiplication. 

    If we fix a certain quaternion of the unit norm $I \in \mathcal{S}$ and consider restriction of $p$ on $\mathbb{R}[I]$, then  behavior of the roots becomes more standard. Since   $\mathbb{R}[I]$ is isomorphic to $\mathbb{C}$, the equality $p_{I}(x) = (x - I)(x + I)$ holds for all $x \in \mathbb{R}[I]$. Hence the set of roots $R(p_{I})$ consists of two elements, $I$ and $-I$.
    
    
    \begin{definition}\label{def poly proj}
        For any given element $I \in \mathcal{S}$ and a polynomial 
        $$P(x) = a_{n}x^n + a_{n-1}x^{n-1}+ \ldots + a_{1}x + a_{0} \in  \mathbb{H}[x],$$
        we denote
        $$P_{I}(x) = b_n x^{n} + b_{n -1} x^{n - 1} + \ldots + b_{1} x + b_{0},$$
        where $b_t$ is an orthogonal projection of $a_{t}$ on  $\mathbb{R}[I]$, $ t = \overline{0, n}$. Similarly, we define
        $$P_{I}^{\perp}(x) = b_n^{\perp}x^{n} + b_{n -1}^{\perp}x^{n - 1} + \ldots + b_{1}^{\perp}x + b_{0}^{\perp},$$
        where $b_{t}^{\perp} = a_{t} - b_{t}$ is the orthogonal complement of $b_{t}$. Let us restrict polynomials $P_{I}$ and $P_{I}^{\perp}$ on the plane $\mathbb{R}[I]$. Accordingly, by the roots of projections we will understand the elements of sets $$R(P_{I}) = \{h \in \mathbb{R}[I] \,\,\, | \,\,\, P_{I}(h) = 0\},$$  
        $$R(P_{I}^{\perp}) = \{h \in \mathbb{R}[I] \,\,\, | \,\,\, P_{I}^{\perp}(h) = 0\}.$$
        
    \end{definition}  
        Note that the nonzero coefficients of $P_{I}^{\perp}$ do not belong to $\mathbb{R}[I]$. However, we can bring it into a good form.
    \begin{proposal}\label{prop perp complex}
         For any polynomial $P_{I}^{\perp}$ there exist $J \in \mathcal{S}$ such that $P_{I}^{\perp} = JQ_{I}$  and   $Q_{I} \in \mathbb{R}[I][x]$.
    \end{proposal}
    Proof:
        Let $P_{I}^{\perp}(x) = b_n^{\perp}x^{n} + b_{n -1}^{\perp}x^{n - 1} + \ldots + b_{1}^{\perp}x + b_{0}^{\perp}$. Since $b^{\perp}_{i}$ lies in an orthogonal complement to $\mathbb{R}[I]$, applying the Lemma \ref{rem 3 I} we obtain the required.
        
    \qed
    
    Let us describe Definition \ref{def poly proj} in more details. The coefficients of the polynomial $P_{I}$ are projections of the coefficients of $P$ on the plane $\mathbb{R}[I]$ with respect to the scalar product \eqref{2}. The coefficients of $P_{I}^{\perp}$ are projections of the coefficients of $P$ onto the orthogonal complement to $\mathbb{R}[I]$. We consider both these  polynomials as functions from  $\mathbb{R}[I]$ to $\mathbb{H}$. Note that for any $h \in \mathbb{R}[I]$ it holds that 
    \begin{equation*} 
        P(h) = P_{I}(h)+P_{I}^{\perp}(h),
    \end{equation*} 
    or in other words 
    $P_{|\mathbb{R}[I]} = P_{I} + P_{I}^{\perp}$. According to Proposition \ref{prop perp complex}, 
    \begin{equation*}
    P_{|\mathbb{R}[I]} = P_{I} + JQ_{I}, 
    \end{equation*}
    where the polynomials $P_{I}$ and $ Q_{I} \in \mathbb{R}[I][x]$, and $J$ is the quaternion of the unit norm introduced in Lemma \ref{rem 3 I}. It is important to underline that the equality $P_{|\mathbb{R}[I]} = P_{I} + JQ_{I}$ is the equality of functions   from $\mathbb{R}[I]$ to~$\mathbb{H}$.
        
    The following proposition describes the relations between the roots of a polynomial and the roots of its projections.

    \begin{proposal}\label{root of slice its intersection of roots projections}
        Let $P \in \mathbb{H}[x]$ and $I \in \mathcal{S}$, then 
        $$R(P_{I}) \cap R(P_{I}^{\perp}) = R(P) \cap \mathbb{R}[I].$$
    \end{proposal}

    Proof:  
    Let $h\in \mathbb{R}[I] \cap R(P)$. By Proposition \ref{prop perp complex} the representation $P(x) = P_{I}(x) + JQ_{I}(x)$ holds, where  $P_{I}, Q_{I}\in \mathbb{R}[I][x]$. Due to associativity, the equality of   quaternions is valid: $P^{\perp}_{I}(q) = JQ_{I}(q)$, for all $q \in \mathbb{R}[I]$. Then we get $P_{I}(h) + JQ_{I}(h) = 0$. According to   Proposition \ref{local complex} quaternions $P(h)$ and $Q(h)$ belongs to $\mathbb{R}[I]$. Then from the Lemma \ref{rem 3 I} it follows  that $P(h) = 0$ and $JQ(h) = 0$. Hence $h$ belongs to the intersection $R(P_{I}) \cap R(P_{I}^{\perp})$.

    On the other hand, if $h \in R(P_{I}) \cap R(P_{I}^{\perp})$, then $P_{I}(h) = 0$ and $P_{I}^{\perp}(h) = 0$. Hence $h$ a root of polynomial $P(x)$.
    
    \qed

To prove the main result  we  consider the roots of the restriction of the polynomial $P$ to $\mathbb{R}[I]$. For this, we need the following corollary.
    
    \begin{consequence}\label{rem from root of slice its intersection of roots projections}
        Let $P \in \mathbb{H}[x]$ and $I \in \mathcal{S}$, then
        $$R(P_{|\mathbb{R}[I]}) = R(P_{I}) \cap R(P_{I}^{\perp}).$$
    \end{consequence}

    Proof:
For any $h \in \mathbb{R}[I]$ 
$$R(P_{|\mathbb{R}[I]})   = \{h \in \mathbb{R}[I] \,\,\, | \,\,\, P(h) = 0\} = R(P) \cap \mathbb{R}[I] .$$
Then the result follows from Proposition \ref{root of slice its intersection of roots projections}.  
        
    \qed
    
    Now, let us understand the structure of the polynomials $P_{I}$. To do this, we need the following lemma.
        
    \begin{lemma}
        For $I \in \mathcal{S}$ the algebras $\mathbb{R}[I][x]$ and $\mathbb{C}[x]$ are isomorphic.
    \end{lemma}
    
    Proof: It follows directly from the Proposition \ref{local complex}.
        
    \qed    

    First, we note that the polynomials of $P_{I}$ lie in the algebra of $\mathbb{R}[I][x]$. And then we will consider the following

    \begin{lemma}
        $\mathbb{R}[I]$ is algebraically closed and any $f \in \mathbb{R}[I][x]$ polynomial of degree at least two satisfies the Gauss-Lucas theorem (Theorem \ref{GL} from this text).
    \end{lemma}
    Proof:
        Since $\mathbb{R}[I]$ and $\mathbb{C}$ are isomorphic and, moreover, isometric, and the algebras of polynomials in these fields are isomorphic, we can apply theorems valid for $\mathbb{C}$ and based on complex numbers of geometric and algebraic properties for $\mathbb{R}[I]$. In particular, the fundamental theorem of algebra and the Gauss-Lucas theorem valid for~$\mathbb{R}[I][x]$.
    \qed
    
    We illustrate the projections of polynomials introduced by the following way.
    \begin{example}
        Let  
        $$P(x) = (1 + \frac{2}{3}i + \frac{1}{3}j + \frac{2}{3}k)x^{3} + (-9 -2i -25j -8k)x^2 +$$
        $$+ (-21 - 115\frac{1}{3}i + 146\frac{1}{3}j -7\frac{1}{3}k)x  + 85 + 367\frac{1}{3}i - 20\frac{1}{3}j + 193\frac{1}{3}k \in \mathbb{H}[x].$$
        Let $g = 1 + i + j$. Then  $g^{-1} = \frac{1 - i - j}{3}$, $I = gig^{-1} = \frac{1}{3}i + \frac{2}{3}j - \frac{2}{3}k \in \mathcal{S}$ and $J = gjg^{-1} = \frac{2}{3}i + \frac{1}{3}j + \frac{2}{3}k$. So the polynomial   $P(x)$ is equal to 
        $$P(x) = (1 + J)x^3 + ((-9 - 12I) + J(-15 - 18I))x^2 + ((-21 + 64I) + J(-33 + 172I))x + $$  $$+(85 - 20I) + J(367 - 194I).$$
        By Lemma \ref{rem 3 I}, we get
        $$P_{I}(x) = x^3 + (-9 - 12I)x^2 + (-21 + 64I)x + 85 - 20I$$
        and, accordingly, 
        $$P_{I}^{\perp}(x) = Jx^{3} + J(-15 -18I)x^{2} + J(-33 + 172I)x + J(367 - 194I).$$
        By the Definition, the roots of these polynomials lie in the plane $\mathbb{R}[I]$. By Definition \ref{def poly proj} the polynomials $P_I$ and $P_I^\perp$ are restricted to the plane $\mathbb{R}[I]$. Subsequently, only the roots lying in this plane are considered. Therefore, the sets of the roots under consideration are equal to $$R(P_{I}) = \{1 + 2I,  3 + 4I,  5 + 6I\}, \,\,\, R(P_{I}^{\perp}) =\{3 + 4I,  5 + 6I, 7 + 8I\}.$$
        
    \end{example}

    This example shows that the results of projections of a polynomial can be nontrivial. Observe that it is not always the case. For example, let us consider a polynomial $P(x) = a_{n}x^{n} + a_{n-1}x^{n-1} + \dots + a_{0}$, where $a_{i} \in \mathbb{R}, \,\,\, i = \overline{0,n}$. Then for all $I \in \mathcal{S}$ we get $P_{I}(x) = a_{n}x^{n} + a_{n-1}x^{n-1} + \dots + a_{0}$ and $P_{I}^{\perp}(x) = 0$. That is, in the case of polynomials with real coefficients, the projections will be trivial.  
    
    \begin{proposal}\label{dsnt matter order}
        For any polynomial $P(x) = \sum\limits_{t = 0}^{n}a_{t}x^{t} \in \mathbb{H}[x]$ and an element $I \in \mathcal{S}$ we have $(P_{I})' = (P')_{I}$.
    \end{proposal}

    Proof: Let $\pi : \mathbb{H} \to \mathbb{R}[I]$ be projection operator on $\mathbb{R}[I]$. Then linearity of $\pi$ implies

        $$(P_{I})' = \sum\limits_{t = 0}^{n} t \pi (a_{t})x^{t - 1} = \sum\limits_{t = 0}^{n} \pi (ta_{t})x^{t - 1} = (P')_{I}$$

    \qed

    The Proposition \ref{dsnt matter order} shows that $P_{I}' = (P_{I})' = (P')_{I}$.

\section{The Gauss-Lucas theorem for quaternions}
        
    \begin{definition} \label{D3.9}
        For a polynomial $P(x) \in \mathbb{H}[x]$ let $$\mathfrak{snm}(P) = \bigcup\limits_{I \in \mathcal{S}} \biggl(\mathbf{conv} (R(P_{I})) \cap \mathbf{conv} (R(P_{I}^{\perp}))\biggr).$$
    \end{definition}

    Let us find the connection between $\mathfrak{sn}(P)$, the snail, proposed by Ghiloni and Perotti and the hull $\mathfrak{snm}(P)$ introduced in   Definition~\ref{D3.9}. To do this, we introduce the notion of cosnail: $\mathfrak{cosn}(P) = \bigcup\limits_{I \in \mathcal{S}} \mathbf{conv} (R(P_{I}^{\perp}))$. Then the following result holds:
    \begin{lemma}  \label{L4.2} For any $P\in \mathbb{H}[x]$ it holds that
      $  \mathfrak{snm}(P) \subseteq \mathfrak{sn}(P) \cap \mathfrak{cosn}(P)  \subseteq  \mathfrak{sn}(P).$ 
    \end{lemma}

    Proof:
    $$\mathfrak{snm}(P) = \bigcup\limits_{I \in \mathcal{S}} \mathbf{conv} (R(P_{I})) \cap \mathbf{conv} (R(P_{I}^{\perp})) \subseteq \bigcup\limits_{I_1, I_2 \in \mathcal{S}} \mathbf{conv} (R(P_{I_1})) \cap \mathbf{conv} (R(P_{I_2}^{\perp})) = $$

    $$ = \bigcup\limits_{I_1\in \mathcal{S}} \biggl( \bigcup\limits_{I_2 \in \mathcal{S}} \Big( \mathbf{conv} (R(P_{I_1})) \cap \mathbf{conv} (R(P_{I_2}^{\perp})) \Big) \biggl) = \bigcup\limits_{I_1\in \mathcal{S}} \biggl( \mathbf{conv} (R(P_{I_1})) \cap \Big (\bigcup\limits_{I_2 \in \mathcal{S}} \mathbf{conv} (R(P_{I_2}^{\perp})) \Big) \biggl) = $$
    $$ = \Big( \bigcup\limits_{I \in \mathcal{S}} \mathbf{conv} (R(P_{I})) \Big) \cap \Big (\bigcup\limits_{I \in \mathcal{S}} \mathbf{conv} (R(P_{I}^{\perp})) \Big) =  \mathfrak{sn}(P) \cap \mathfrak{cosn}(P)  \subseteq  \mathfrak{sn}(P).$$
Here we use the property $(\bigcup\limits_{\alpha}A_{\alpha}) \cap B = \bigcup\limits_{\alpha}(A_\alpha \cap B)$ (see, \cite[\S 1.2]{A}).
    
    \qed

    We get that $\mathfrak{snm}(P)$ lies in the intersection of two sets. The first of which is $\bigcup\limits_{I \in \mathcal{S}} \mathbf{conv} (R(P_{I}))$ the snail which is already familiar to us. The second set is $\bigcup\limits_{I \in \mathcal{S}} \mathbf{conv} (R(P_{I}^{\perp}))$. This set differs from the snail since  instead of orthogonal projections on the plane $\mathbb{R}[I]$ we take orthogonal projections on orthogonal complement to $\mathbb{R}[I]$. It is called a cosnail in order to  emphasize this duality.  
    Thus, for any polynomial $P$ it holds   $\mathfrak{snm} (P) \subseteq \mathfrak{sn} (P) \cap \mathfrak{cosn} (P)$. In particular, this implies the inclusion $\mathfrak{snm}(P) \subseteq \mathfrak{sn}(P)$. In Corollary \ref{strictly enabling} we will show that this inclusion can be strict. We also note that $\mathfrak{snm} (P)$ does not differ much from $\mathfrak{sn} (P) \cap \mathfrak{cosn}(P)$ as the following proposition shows.

    \begin{proposal} For any polynomial $P$ we have  $\Bigl( \mathfrak{sn} (P) \cap \mathfrak{cosn} (P) \Bigl) \setminus \mathfrak{snm} (P) \subseteq \mathbb{R}$.
    \end{proposal}
    Proof:
By   Lemma~\ref{L4.2} 
  $$\mathfrak{snm}(P)  \subseteq   \mathfrak{sn} (P) \cap \mathfrak{cosn} (P) = \bigcup\limits_{I_1, I_2 \in \mathcal{S}} \mathbf{conv} (R(P_{I_1})) \cap \mathbf{conv} (R(P_{I_2}^{\perp})). $$

In the case  $I_{1}=\pm I_{2}$, we get by Definition \ref{D3.9}  that $\mathbf{conv} (R(P_{I_1})) \cap \mathbf{conv} (R(P_{I_1}^{\perp})) \subseteq \mathfrak{snm}(P)$.

In the case  $I_1 \neq \pm I_2$,   the intersection $\mathbf{conv} (R(P_{I_1})) \cap \mathbf{conv} (R(P_{I_2}^{\perp})) \subseteq \mathbb{R}[I_1] \cap \mathbb{R}[I_2] = \mathbb{R}$, since $ R(P_{I}) \subseteq \mathbb{R}[I]$ for any $I$ and also $\mathbf{conv} (R(P_{I})) \subseteq \mathbb{R}[I]$.   

These two cases together  imply that $\bigl( \mathfrak{sn} (P) \cap \mathfrak{cosn} (P) \bigr) \setminus \mathfrak{snm} (P) \subseteq \mathbb{R} $. 
        
    \qed
    
    Now,  using Definition \ref{D3.9} we can reformulate  Theorem \ref{SQGL} in the following way and proceed to prove it.    
    \begin{theorem}\label{ssn}
        For any polynomial $P \in \mathbb{H}[x]$
        \begin{equation*}
            R(P') \subseteq \mathfrak{snm}(P).   
        \end{equation*}
        
    \end{theorem}

    Proof:
    We fix $I \in \mathcal{S}$ and consider the decomposition $P_{|\mathbb{R}[I]} = P_{I} + P_{I}^{\perp}$ provided by Definition \ref{def poly proj}. Then Proposition \ref{prop perp complex} implies that   $P_{|\mathbb{R}[I]} = P_{I} + JQ_{I}$, where $J \in \mathcal{S}$ and $Q_{I}$ has coefficients from $\mathbb{R}[I]$.
        
    Hence $P'_{|\mathbb{R}[I]} = P_{I}' + JQ_{I}'$. From Corollary \ref{rem from root of slice its intersection of roots projections} we have $a \in R(P'_{|\mathbb{R}[I]})$ if and only if $a \in R(P'_{I}) \cap R({P'}_{I}^{\perp})$, or equivalently $a \in R(P'_{I}) \cap R(Q'_{I})$. By applying Theorem \ref{GL} to the polynomials $P_{I}$ and $Q_{I}$, we get

        \begin{equation*}
            R(P_{I}') \subseteq \mathbf{conv} (R(P_{I})), 
        \end{equation*}

        \begin{equation*}
            R(Q_{I}') \subseteq \mathbf{conv} (R(Q_{I})) =\mathbf{conv} (R(P_{I}^{\perp})),
        \end{equation*}

        \begin{equation*}
            R((P')_{| \mathbb{R}[I]}) = R(P_{I}') \cap R(Q_{I}') \subseteq \mathbf{conv} (R(P_{I})) \cap \mathbf{conv} (R(P_{I}^{\perp})).
        \end{equation*}
        Each root lies in a certain $\mathbb{R}[I]$, hence
            
        $$R(P') = \underset{I \in \mathcal{S}}{\bigcup}R((P')_{| \mathbb{R}[I]}) \subseteq \bigcup\limits_{I \in \mathcal{S}} \mathbf{conv} (R(P_{I})) \cap \mathbf{conv} (R(P_{I}^{\perp})).$$
        
    \qed

By Lemma~\ref{L4.2} $\mathfrak{snm}(P) \subseteq \mathfrak{sn}(P)$. Let us show that the inclusion can be strict.  

    \begin{example}\label{Ex:strictly enabling}
 For the polynomial $P(x) = x^3 - (i + j + k) x^2 + (k + i - j)x + 1 \in \mathbb{H}[x]$        the inclusion $\mathfrak{snm}(P) \subset \mathfrak{sn}(P)$ is strict.
    \end{example}

    Proof:
        Fix the polynomial $P(x) = x^3 - (i + j + k) x^2 + (k + i - j)x + 1 \in \mathbb{H}[x]$. Since $i \in \mathcal{S}$, we can consider the plane $\mathbb{R}[i]$. Consider the projection $P_i$ of $P(x)$ on $\mathbb{R}[i]$ and its orthogonal complement $P_{i}^{\perp}(x)$: 
        $$P_{i}(x) = x^3 - ix^2 + ix +1$$
        and
        $$P_{i}^{\perp}(x) = -(j + k)x^2 - (j - k)x.$$
        Decompose both polynomials into the products of irreducible complex polynomials, irreducible polynomials over the plane $\mathbb{R}[i]$:
        $$P_{i}(x) = (x - (\cos{\frac{\pi}{4}} + i \sin{\frac{\pi}{4}}))(x - (\cos{\frac{5\pi}{4}} + i \sin{\frac{5\pi}{4}}))(x - i)$$
        and
        $$P_{i}^{\perp}(x) = -j((1 - i)x + (1 + i))x.$$
        On the one hand, this means that set of roots of polynomial $P_{i}(x)$ is equal to

        $$R((P_{i})_{| \mathbb{R}[i]}) = \{i, (\cos{\frac{\pi}{4}} + i \sin{\frac{\pi}{4}}), (\cos{\frac{5\pi}{4}} + i \sin{\frac{5\pi}{4}})\}.$$
        Thus, $\mathbf{conv} (R(P_{i}))$ is a triangle with vertices at the listed points on the complex plain. On the other hand,

        $$R(P_{i}^{\perp}) = \{0, -i\}.$$
        Then $\mathbf{conv} (R(P_{i}^{\perp}))$ is a segment with its ends at the points $0$ and $-i$. Since the triangle $\mathbf{conv} (R(P_{i}))$ intersects the imaginary axis along the segment that ends at $i$ and $0$,  we get

        $$ \{0\} = \mathbf{conv} (R(P_{i})) \cap \mathbf{conv} (R(P_{i}^{\perp})) \subsetneq  \mathbf{conv} (R(P_{i})).$$
        Set $\mathbf{conv} (R(P_{i}))$ belong to $\mathfrak{sn}(P)$. Since $\mathbf{conv} (R(P_{i})) \not\subseteq \mathbb{R}[I]$ for $I \neq \pm i$, then $\mathbf{conv} (R(P_{i})) \not\subseteq \mathbf{conv} (R(P_{I})) \cap \mathbf{conv} (R(P_{I}^{\perp}))$ for $I \neq \pm i$. At the same time $\mathbf{conv} (R(P_{i})) \not\subseteq \mathbf{conv} (R(P_{i})) \cap \mathbf{conv} (R(P_{i}^{\perp}))$. This implies that for the polynomial $P(x) = x^3 - (i + j + k) x^2 + (k + i - j)x + 1$ the inclusion is strict: 
        
        \begin{equation*}\label{eq 5}
            \mathfrak{snm}(P) = \bigcup\limits_{I \in \mathcal{S}} \mathbf{conv} (R(P_{I})) \cap \mathbf{conv} (R(P_{I}^{\perp})) \subsetneq \mathfrak{sn}(P).
        \end{equation*}

        \begin{figure}[h]
            \centering
            \includegraphics[width=0.5\textwidth]{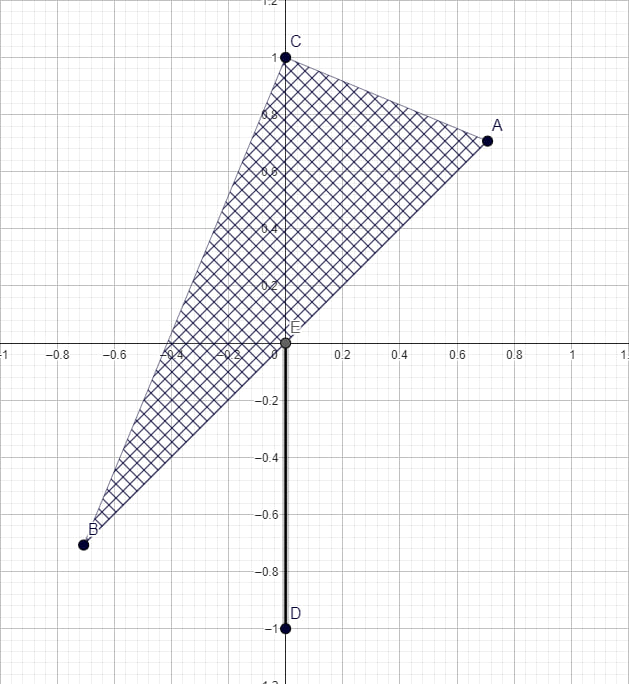}
            \caption{\scriptsize Set $\mathbf{conv} (R(P_{i}))$ is a triangle $ABC$. 
            Set $\mathbf{conv} (R(P_{i}^{\perp}))$ is a segment $DE$. \normalsize}
        \end{figure}
  
    \qed

    \begin{consequence}
    \label{strictly enabling}
        Inclusion $\mathfrak{snm}(P) \subseteq \mathfrak{sn}(P)$ can be strict.
\end{consequence}
\begin{proof}
    Follows from Example~\ref{Ex:strictly enabling}.
\end{proof}

    \begin{remark}
        The set $\mathfrak{snm}(P)$ from the proof  of Example~\ref{Ex:strictly enabling} is not empty. Indeed, for the element $j \in \mathcal{S}$ the projections $P_{j}$ and $P_{j}^{\perp}$ of the polynomial $P$ are
        $$P_{j}(x) = x^{3} -jx^{2} -jx + 1$$
        and
        $$P_{j}^{\perp}(x) = -(i + k)x^2 + (i + k)x.$$
        Then we get
        $$\{-1\} \in R(P_{j}) \cap R(P_{j}^{\perp}).$$
    \end{remark}

    The convex hull we have introduced is superior not only because it provides a thinner and more precise geometric structure, as demonstrated by Lemma \ref{L4.2} and Proposition \ref{strictly enabling}. The key point is that when we consider only the orthogonal projections of a polynomial, we lose approximately half of the information about it. For example, if we consider the projection of the polynomial $j(x^2+1)$ onto the plane $\mathbb{R}[i]$, its projection turns out to be zero, and we lose all information about its roots. In contrast, Proposition \ref{root of slice its intersection of roots projections}  shows that by considering both orthogonal projections and orthogonal complements, we can recover the roots of the original polynomial. Thus, we obtain a more symmetric and more natural geometric object. 

    Thus, our version of the Gauss–Lucas theorem is more refined and relies more heavily on the geometric properties of the polynomial.
    
\section{Estimation of the norm of the roots}
Now we are ready to discuss the applications of the above results. For the polynomial $Q(x) = \sum\limits_{t = 0}^{n} a_{t}x^{t} \in \mathbb{C}[x]$ the following bound is known,~\cite[\S 8.1]{QRGS}:

    \begin{equation}\label{mod est}
        \max\limits_{z \in R(Q)} |z| \leq  |a_{n}|^{-1} \sqrt{\sum\limits_{t = 0}^{n} |a_{t}|^{2}}.
    \end{equation}
    
    We will introduce the norm of the polynomial, which will be used to estimate the norm of the roots.
    
    \begin{definition}
        For the polynomial $P(x)= \sum\limits_{t = 0}^{n} a_{t}x^{t} \in \mathbb{H}[x]$ consider the function $C: \mathbb{H}[x] \to \mathbb{R} \cup \{+\infty \}$  defined as follows:
        \begin{equation*}
            С(P)=
            \begin{cases}
                +\infty, & \text{if}\ P=const. \\
               |a_{n}|^{-1} \sqrt{\sum\limits_{t = 0}^{n} |a_{t}|^{2}}, & \text{if}\ n \geq 1\ \text{and}\ a_{n} \neq 0.
            \end{cases}
        \end{equation*}
    \end{definition}
 It is straightforward that the function $C$ is a  norm on the polynomial space. 
    In \cite{GRPA}, using the inequality (\ref{mod est}), the following estimate for the roots of the derivative was proved.

    \begin{lemma}\cite[Proposition 6]{GRPA}
        For any polynomial $P(x)= \sum\limits_{t = 0}^{n} a_{t}x^{t} \in \mathbb{H}[x]$ such that $n \geq 1$ and $a_{n} \neq 0$ we have
        \begin{equation}\label{canon estimate module dervative roots}
            \max\limits_{x \in R(P')} |x| \leq \sup\limits_{I \in \mathcal{S}} C(P_{I}).
        \end{equation}
        
    \end{lemma}

    By Theorem \ref{ssn} this estimate can be strengthened.
    
    \begin{lemma}
        For any polynomial $P(x)= \sum\limits_{t = 0}^{n} a_{t}x^{t} \in \mathbb{H}[x]$ such that $n \geq 1$ and $a_{n} \neq 0$ the following inequality is true:
        \begin{equation}\label{estimate module dervative roots}
            \max\limits_{x \in R(P')} |x| \leq \sup\limits_{I \in \mathcal{S}} \{\min (C(P_{I}), C(P_{I}^{\perp})\}.
        \end{equation}
    \end{lemma}

    Proof:
        It follows from  Theorem \ref{ssn} that the root of $h\in R(P')$ belongs to a certain set of the form $\mathbf{conv} (R(P_{I})) \cap \mathbf{conv} (R(P_{I}^{\perp}))$. Hence $h$ belong to $\mathbf{conv} (R(P_{I}))$, and $\mathbf{conv} (R(P_{I}^{\perp}))$ simultaneously. Due to the estimate (\ref{mod est}) we get $|h| \leq C(P_{I})$ and $|h| \leq C(P_{I}^{\perp})$,  the statement of the lemma follows.
        
    \qed

   Now we show that our estimate indeed improve the statement of \cite[Proposition 6]{GRPA}. To do this, consider

    \begin{proposal}
        The inequality $\sup\limits_{I \in \mathcal{S}} \{\min (C(P_{I}), C(P_{I}^{\perp})\} \leq \sup\limits_{I \in \mathcal{S}} C(P_{I})$ can be strict.
    \end{proposal}
    
    Proof: 
        Consider the polynomial $P(x) = x^3 + ix^2 + 3jx$. For any $I \in \mathcal{S}$ the projections of the polynomial $P(x)$ are:
        $$P_{I}(x) = x^3 + (i,I) \cdot I x^2 + 3(j, I) \cdot Ix,$$
        $$P_{I}^{\perp}(x) = (i -(i,I) \cdot I) x^2 + 3(j - (j, I) \cdot I)x.$$

        If we substitute $I = j$, we get $C(P_I) = \sqrt{1 + |(i,j)|^2 + 9 |(j,j)|^2} = \sqrt{10}$. Hence $ \sqrt{10} \leq \sup\limits_{I \in \mathcal{S}} C(P_{I})$.

        Let us bound $\sup\limits_{I \in \mathcal{S}} \{\min (C(P_{I}), C(P_{I}^{\perp})\}$ from the above. 
        
        If $I = i$, then $C(P_I) = \sqrt{2}$ and $\min (C(P_{I}), C(P_{I}^{\perp}))$ do not exceed $\sqrt{2}$. Next, let $I = \alpha i + \beta j + \gamma k$ belong to $\mathcal{S}$ and not equal to $i$. Then the norms of projections of the polynomial $P(x)$ are equal to:
        $$С(P_{I}) = \sqrt{1 + |(i,I)|^2 + 9|(j, I)|^2} = \sqrt{1 + \alpha^2 + 9 \beta^2},$$
        $$C(P_{I}^{\perp}) = |i - (i,I) \cdot I|^{-1} \cdot \sqrt{|i - (i, I) \cdot I|^2 + 9|j - (j, I) \cdot I|^2 } = \sqrt{1 + 9 \cdot \frac{1 - \beta^2}{1 - \alpha^2}}.$$
        Since $I \in \mathcal{S}$, then $\alpha^2 + \beta^2 \leq 1$. General situation splits into the following subcases.
        \begin{enumerate}
            \item If $\alpha \geq \frac{1}{2}$, then $C(P_{I}) \leq \sqrt{1 + \alpha^2 + 9(1 -\alpha^2)} \leq \sqrt{8}.$
            \item If $\alpha \leq \frac{1}{2}$, then 
            $$C(P_{I}) \leq \sqrt{1 + \frac{1}{4} + 9\beta^2},$$
            $$C(P_{I}^{\perp}) \leq \sqrt{1 + 9 \cdot \frac{1 - \beta^2}{1 - \frac{1}{4}}} = \sqrt{13 - 12 \beta^2}.$$

            \begin{enumerate}
                \item If $\beta \leq \frac{1}{2}$, then $C(P_{I}) \leq \sqrt{1 + \frac{10}{4}}.$
                \item If $\beta \geq \frac{1}{2}$, then $C(P_{I}^{\perp}) \leq \sqrt{9}.$
            \end{enumerate}
        \end{enumerate}

        In total, we get: 
        $$\sup\limits_{I \in \mathcal{S}} \{\min (C(P_{I}), C(P_{I}^{\perp})\} \leq \sqrt{9} < \sqrt{10} \leq \sup\limits_{I \in \mathcal{S}} C(P_{I}).$$

        That is, for the polynomial $P(x) = x^3 + ix^2 + 3jx$, the estimate (\ref{estimate module dervative roots}) is less than or equal to $\sqrt{9}$, and the estimate (\ref{canon estimate module dervative roots})   is greater than $\sqrt{10}$. Therefore, for the polynomial $P(x)$, the inequality (\ref{estimate module dervative roots}) estimates the absolute value of the roots of $P'(x)$ more precisely than (\ref{canon estimate module dervative roots}).
        
    \qed

\end{document}